\documentclass[12pt]{article}
\usepackage{amsfonts}

\usepackage{latexsym}
\usepackage{amssymb}
\usepackage{amsmath}
\usepackage{amsthm}
\usepackage[mathscr]{eucal}
\usepackage{graphicx}
\usepackage{hyperref}
\usepackage{caption}
\usepackage{subcaption}

\newtheorem{Lemma}{Lemma}

\newtheorem{Theorem}[Lemma]{Theorem}

\renewcommand{\qed}{\hfill{\ \ \rule{2mm}{2mm}} \vspace{0.2in}}

\begin{document}

\title{Deviation estimates for Eulerian edit numbers of random graphs}
\author{ \textbf{Ghurumuruhan Ganesan}
\thanks{E-Mail: \texttt{gganesan82@gmail.com} } \\
\ \\
Institute of Mathematical Sciences, HBNI, Chennai}
\date{}
\maketitle

\begin{abstract}
Consider the random graph~\(G(n,p)\) obtained by allowing each edge in the complete graph on~\(n\) vertices to be present with probability~\(p\) independent of the other edges. In this paper, we study the minimum number of edge edit operations needed to convert~\(G(n,p)\) into an Eulerian graph. We obtain deviation estimates for three types Eulerian edit numbers based on whether we perform only edge additions or only edge deletions or a combination of both and show that with high probability, roughly~\(\frac{n}{4}\) operations suffice in all three cases.

\vspace{0.1in} \noindent \textbf{Key words:} Eulerian edit numbers, random graphs, deviation estimates.

\vspace{0.1in} \noindent \textbf{AMS 2000 Subject Classification:} Primary: 60J10, 60K35.
\end{abstract}

\bigskip

\renewcommand{\theequation}{\thesection.\arabic{equation}}
\setcounter{equation}{0}
\section{Introduction} \label{intro}
The problem of converting a given graph~\(G\) into an Eulerian graph is of great theoretical and practical importance. For example, Boesch et al. (1977) show that~\(G\) can be extended into Eulerian graph if and only if~\(G\) is not an odd complete bipartite graph and also determine the minimum number of edge additions needed. Random graphs have frequently been used to estimate edit distance for hereditary properties (i.e. properties that are invariant under removal of vertices) of deterministic graphs. For example Alon and Stav (2008) show that a random graph essentially achieves the maximum distance for any hereditary property and determine the ``critical value" of the edge probability for  this to occur. For more details, we also refer to the  recent survey paper of Martin (2016).

In this paper, we study the question of Eulerian extendability/reducability for random graphs based on the Erd\H{o}s-R\'enyi model. In our main deviation estimate result in Theorem~\ref{thm_eul_edit} stated below, we show that with high probability, i.e. with probability converging to one as~\(n \rightarrow \infty,\) roughly~\(\frac{n}{4}\) edge edit operations are enough to convert a random graph into an Eulerian graph.

\subsection*{Eulerian edit numbers} \label{eul_edit}
Let~\(G = (V,E)\) be a graph with vertex set~\(V\) and edge set~\(E.\) An edge between vertices~\(u\) and~\(v\) is denoted as~\((u,v)\) and in this case~\(u\) and~\(v\) are said to be adjacent to each other. A sequence of vertices~\({\cal W} := (u_1,u_2,\ldots,u_t), t \geq 3\) is said to be a \emph{circuit} if~\(u_i\) is adjacent to~\(u_{i+1}\) for each~\(1 \leq i \leq t-1\) and~\(u_t\) is adjacent to~\(u_1.\) We say that~\({\cal W}\) is an Eulerian circuit if each edge of the graph~\(G\) occurs exactly once in~\({\cal W}.\) The graph~\(G\) is said to be an \emph{Eulerian} graph if~\(G\) contains an Eulerian circuit. By definition any Eulerian graph is connected and contains at least three edges.

A graph~\(G(1)\) obtained by either an addition or a removal of a single edge in~\(G\) is defined to be a~\(1-\)\emph{edit graph} of~\(G.\) Iteratively, for~\(i \geq 1\) if~\(G(i)\) is an~\(i-\)edit graph of~\(G,\) then any~\(1-\)edit of~\(G(i)\) is defined to be an~\((i+1)-\)edit graph of~\(G.\) For completeness the~\(0-\)edit graph of~\(G\) is the graph~\(G\) itself. We define the \emph{Eulerian edit number}~\(N_{edit}(G)\) of~\(G\) to be the smallest integer~\(k \geq 0\) such that some~\(k-\)edit graph~\(G(k)\) of~\(G\) is Eulerian. The number~\(N_{edit}(G)\) is well-defined since if~\(G\) contains~\(n\) vertices and~\(m\) edges, then~\(N_{edit}(G) \leq m+n\) as we can remove all the edges of~\(G\) and then add a cycle containing all the~\(n\) vertices, to obtain an Eulerian graph.

We say that~\(G\) is \emph{Eulerian extendable} if~\(G\) can be converted into an Eulerian graph by adding a finite number of edges. We define the \emph{Eulerian extension number}~\(N_{ext}(G)\) of an Eulerian extendable graph~\(G\) to be the minimum number of edges that must be added to~\(G\) so that the resulting graph is Eulerian. We say that~\(G\) is \emph{Eulerian reducable} if~\(G\) can be converted into an Eulerian graph containing at least one edge, by removing a finite number of edges. We define the \emph{Eulerian reduction number}~\(N_{red}(G)\) of an Eulerian reducable graph~\(G\) to be the minimum number of edges that must be removed from~\(G\) so that the resulting graph is Eulerian.

Let~\(K_n\) be the complete graph on~\(n\) vertices and let~\(\{X_e\}_{e \in K_n}\) be independent and identically distributed (i.i.d.) random variables indexed by the edge set of~\(K_n,\) defined on the probability space~\((\Omega,{\cal F},\mathbb{P})\) and having distribution
\[\mathbb{P}(X_e = 1) = p = 1-\mathbb{P}(X_e = 0).\] Let~\(H = G(n,p)\) be the random graph obtained by allowing edge~\(e \in K_n\) to be present if~\(X_e =1\) and absent otherwise. We have the following result regarding the previously defined Eulerian numbers of~\(H.\) All constants throughout do not depend on~\(n.\)
\begin{Theorem}\label{thm_eul_edit}
Suppose
\begin{equation}\label{p_cond2}
\frac{(\log{n})^2}{\sqrt{n}}\leq p \leq 1-\frac{(\log{n})^2}{\sqrt{n}}
\end{equation}
and for~\(\eta > 0\) let~\(E_{edit} = E_{edit}(\eta)\) be the event that
all three Eulerian numbers~\(N_{edit}(H),N_{ext}(H)\) and~\(N_{red}(H)\) lie between~\(\frac{n}{4} - n^{\frac{1}{2}+\eta} \) and~\(\frac{n}{4} + n^{\frac{1}{2}+\eta}.\)

For every~\(r >0\) there exists an integer~\(n_0 = n_0(\eta,r) \geq 1\) such that if~\(n \geq n_0\) then
\begin{equation}\label{eul_edit_est}
\mathbb{P}(E_{edit}) \geq 1-\frac{1}{n^{r}}.
\end{equation}
\end{Theorem}
With high probability therefore, the random graph~\(H\) can be converted into an Eulerian graph by adding or removing roughly~\(\frac{n}{4}\) edges. For comparison, we get from~(\ref{p_cond2}) that the expected number of edges~\(m := {n \choose 2} \cdot p\) in the graph~\(H\) is at least~\(n^{\frac{3}{2}} \cdot (\log{n}).\) Consequently, from standard Binomial deviation estimates (Corollary~\(A.1.14,\) pp. 312, Alon and Spencer (2008)) we know that with high probability~\(H\) contains at least~\(n^{\frac{3}{2}}\) edges.

The paper is organized as follows: In Section~\ref{prelim}, we collect the preliminary results regarding the graph~\(H\) used in the proof of Theorem~\ref{thm_eul_edit} and in Section~\ref{thm_proof}, we prove Theorem~\ref{thm_eul_edit}.


\renewcommand{\theequation}{\thesection.\arabic{equation}}
\setcounter{equation}{0}
\section{Preliminaries} \label{prelim}
We begin with a result that estimates the number of odd degree vertices in~\(H = G(n,p).\) Let~\(T\) denote the number of odd degree vertices in~\(H\) and let~\(\mu = n \cdot \nu = \mathbb{E}T\) be the mean of~\(T.\)
\begin{Lemma}\label{vert_lem}
Suppose
\begin{equation}\label{p_cond}
\frac{(\log{n})^2}{n}\leq p \leq 1-\frac{(\log{n})^2}{n}.
\end{equation}
For every even integer~\(q \geq 2\) there are constants~\(D = D(q) > 0\) and\\\(n_0  =n_0(q)\) such that if~\(n \geq n_0(q)\) then
\begin{equation}\label{mu_est}
\frac{n}{2} - \frac{1}{n^{q}} \leq \mu \leq \frac{n}{2} + \frac{1}{n^{q}}\;\;\;\text{   and   }\;\;\;\mathbb{E}(T-\mu)^{q} \leq D \cdot n^{\frac{q}{2}}.
\end{equation}
\end{Lemma}
\emph{Proof of Lemma~\ref{vert_lem}}: For a constant integer~\(b \geq 1\) we first define\\\(\epsilon_b :=\frac{1}{2}(1-2p)^{n-b}\) and show that if~(\ref{p_cond}) holds, then
\begin{equation}\label{ep_est}
|\epsilon_b| \leq \frac{1}{n^{q+1}}
\end{equation}
for all~\(n\) large. Indeed if~\(p \leq \frac{1}{2},\) then using the lower bound~\(p\geq \frac{ (\log{n})^2}{n}\) in~(\ref{p_cond}) we have that
\[|\epsilon_b| = (1-2p)^{n-b} \leq e^{-2(n-b)p} \leq \frac{1}{n^{q+1}}\] for all~\(n\) large.  On the other hand, if~\(p > \frac{1}{2},\) then using the upper bound\\\(p \leq 1- \frac{(\log{n})^2}{n}\) in~(\ref{p_cond}) we get
\[|\epsilon_b| = (2p-1)^{n-b} = (1-2(1-p))^{n-b} \leq e^{-2(n-b)(1-p)} \leq \frac{1}{n^{q+1}}\] for all~\(n\) large.
This proves~(\ref{ep_est}).

Next using~(\ref{ep_est}), we prove the estimate for the mean in~(\ref{mu_est}). Denoting~\(T_i\) to be the indicator function of the event that vertex~\(i\) has odd degree we get
\begin{equation}
\mathbb{P}(T_i = 1) = \sum_{\stackrel{0 \leq k \leq n-1}{k \text{ odd }}} {n-1 \choose k} p^{k}(1-p)^{n-1-k}  = \frac{1}{2} - \epsilon_1, \label{fj_comp}
\end{equation}
since the middle term in~(\ref{fj_comp}) equals~\(\frac{1}{2}\left((1-p+p)^{n-1} - ((1-p)-p)^{n-1}\right).\)
From~(\ref{fj_comp}),~(\ref{ep_est}) and the fact that~\(\mathbb{E}T = \sum_{i=1}^{n}\mathbb{E}T_i = \sum_{i=1}^{n} \mathbb{P}(T_i = 1),\) we get the first estimate in~(\ref{mu_est}). For future use we also get from~(\ref{fj_comp}) that for any~\(b \geq 1\) the product
\begin{equation}\label{tj_ind}
\prod_{i=1}^{b} \mathbb{P}(T_i=1) = \left(\frac{1}{2} - \epsilon_1\right)^{b} = \frac{1}{2^{b}} + U_b,
\end{equation}
where~\(|U_b| \leq  \sum_{k=1}^{b} {b \choose k} |\epsilon_1|^k \leq 2^{b} |\epsilon_1| \leq \frac{D_1}{n^{q+1}}\) for all~\(n\) large and some constant~\(D_1 > 0,\) by~(\ref{ep_est}).

We prove the second estimate in~(\ref{mu_est}) in two steps. In the first step, we let~\(b \geq 1\) be an integer constant and use~(\ref{ep_est}) to show that any~\(b\) random variables from~\(\{T_i\}_{1 \leq i \leq n}\) are nearly independent. We then expand~\(|T-\mu|^{q}\) and use the near-independence estimate term-by-term to obtain the second bound in~(\ref{mu_est}). \\
\underline{\emph{Step 1}}: Let~\(H_b\) be the induced subgraph of the random graph~\(H\) formed by the vertices~\(\{1,2,\ldots,b\}\) and write
\begin{equation}\label{tih}
\bigcap_{i=1}^{b} \{T_i = 1\} = \bigcup_{\Gamma} \{H_b = \Gamma\} \bigcap \bigcap_{i=1}^{b} E_i(\Gamma)
\end{equation}
where~\(E_i(\Gamma)\) is the event that the number of neighbours of~\(i\) among the vertices~\(\{b+1,\ldots,n\}\) in the graph~\(H\) has the opposite parity of the number of neighbours of~\(i\) in the graph~\(\Gamma\) i.e. if the vertex~\(i\) has even number of neighbours in~\(\Gamma\) then the number of neighbours of~\(i\) in the set~\(\{b+1,\ldots,n\}\) is odd and vice versa.


For any graph~\(\Gamma,\) the events~\(\{E_i(\Gamma)\}\) are conditionally independent given\\\(H_b = \Gamma\) and so we have from~(\ref{fj_comp}) that
\begin{equation}\label{ei_gam}
\mathbb{P}\left(\bigcap_{i=1}^{b} E_i(\Gamma) \mid H_b = \Gamma\right) = \prod_{i=1}^{b}\left(\frac{1}{2} + q_i(\Gamma)\right)
\end{equation}
where~\(q_i(\Gamma) =\epsilon_b\) if vertex~\(i\) has odd degree in~\(\Gamma\) and~\(q_i(\Gamma) = -\epsilon_b\) otherwise. Expanding the product in~(\ref{ei_gam}) and using~(\ref{ep_est}), we get
\begin{equation}\label{ei_gam2}
\mathbb{P}\left(\bigcap_{i=1}^{b} E_i(\Gamma) \mid H_b = \Gamma\right) = \frac{1}{2^{q}} + r_i(\Gamma)
\end{equation}
where~\(|r_i(\Gamma)| \leq \frac{D_2}{n^{q+1}}\) for all~\(n\) large and some constant~\(D_2 > 0.\) Averaging over~\(\Gamma\) and using~(\ref{tih}) we therefore have
\begin{equation}\label{ptq_est}
\mathbb{P}\left(\bigcap_{i=1}^{b} \{T_i = 1\}\right)  = \frac{1}{2^{b}} + R_b,
\end{equation}
where~\(R_b := \sum_{\Gamma} r_i(\Gamma) \mathbb{P}(H_b = \Gamma)\) is bounded above by~\(|R_b| \leq  \frac{D_2}{n^{q+1}}\) for all~\(n\) large, by the estimate for~\(r_i(\Gamma)\) in~(\ref{ei_gam2}). From~(\ref{ptq_est}) and~(\ref{tj_ind}) we therefore get that
\begin{equation}\label{et_est2}
\left| \mathbb{P}\left(\bigcap_{i=1}^{b} \{T_i = 1\}\right) - \prod_{i=1}^{b}\mathbb{P}\left(T_i =1\right)\right| \leq \frac{D_3}{n^{q+1}}
\end{equation}
for some constant~\(D_3 > 0.\) In other words, for any constant integer~\(b \geq 1\) the random variables~\(\{T_i\}_{1 \leq i \leq b}\) are nearly independent in the sense of~(\ref{et_est2}).

\emph{\underline{Step 2}}: We now use~(\ref{et_est2}) to obtain the second estimate in~(\ref{mu_est}). Recalling that~\(T = \sum_{i=1}^{n}T_i\) with mean~\(\mu = \mathbb{E}T =n \cdot \nu,\) we write
\begin{equation}\label{et_mu}
\mathbb{E}(T-\mu)^{q} = \sum_{i_1,\ldots,i_q} \mathbb{E}\prod_{j=1}^{q} \left(T_{i_j}-\nu\right),
\end{equation}
where every term in the summation in~(\ref{et_mu}) is bounded above by~\(1\) since
by the AM-GM inequality we have
\begin{equation}\label{first_est}
\mathbb{E}\prod_{j=1}^{q} |T_{i_j}-\nu| \leq \frac{1}{q} \sum_{j=1}^{q} \mathbb{E}|T_{i_j}-\nu|^{q} = \mathbb{E}|T_1-\nu|^q \leq 1.
\end{equation}

For a general term in~(\ref{et_mu}), the multiset~\(\{i_1,\ldots,i_q\}\) may have repeated entries and so we write~\(\{i_1,\ldots,i_q\} := \{b_1\cdot a_1,\ldots,b_l \cdot a_l\}\) to denote that the entry~\(a_j\) appears~\(b_j\) times in~\(\{i_1,\ldots,i_q\}.\) Defining~\({\cal B}(i_1,\ldots,i_q) := \min_{1 \leq j \leq l} b_l,\) we let~\({\cal I}_1\) be the set of all multisets~\(\{i_1,\ldots,i_q\}\) for which~\({\cal B}(i_1,\ldots,i_q) = 1\) and let~\({\cal I}_2\) be the remaining set of~\(q-\)element multisets where each element occurs at least twice and split the summation in~(\ref{et_mu}) into two terms~\(I_1\) and~\(I_2\)
where
\begin{equation}\label{etq_split}
I_k := \sum_{\{i_1,\ldots,i_q\} \in {\cal I}_k} \mathbb{E}\prod_{j=1}^{q} \left(T_{i_j}-\nu\right)
\end{equation}
for~\(k=1,2.\) If~\({\cal B}(i_1,\ldots,i_q) \geq 2\) then there are at most~\(\frac{q}{2}\) distinct terms in~\(\{i_1,\ldots,i_q\}\) and the number of choices for such~\(\{i_1,\ldots,i_q\}\) is at most~\({n \choose \frac{q}{2}} \leq n^{\frac{q}{2}}.\) Using the upper bound~(\ref{first_est}) we therefore have that~
\begin{equation}\label{i2_est}
|I_2| \leq n^{\frac{q}{2}}.
\end{equation}

To evaluate the terms in~\({\cal I}_1,\) we write \[\prod_{j=1}^{q} \left(T_{i_j}-\nu\right) = \prod_{j=1}^{l}\left(T_{a_j}-\nu\right)^{b_j} := Q_l\] and show that if~\({\cal B}(i_1,\ldots,i_q)=1\) then the expected value of~\(\mathbb{E}Q_l\) is very small. This is intuitive since if the~\(T_i's\) were independent then in fact we would have that~\(\mathbb{E}Q_l = 0.\) Indeed suppose the element~\(1\) occurs exactly once in~\(\{i_1,\ldots,i_q\}\) so that \(Q_l = (T_1-\nu) \cdot X\) where~\(X = \prod_{j=1}^{q-1} \left(T_{i_j}-\nu\right).\)  Because the~\(T_i's\) are indicator functions, any term~\(w\) in the expansion of~\(X\) is of the form~\(w = (-\nu)^{a} \cdot \prod_{k=1}^{b} T_{l_k}\) for integer~\(0 \leq a \leq q-1\) and~\(0 \leq b \leq q-1\) distinct integers~\(l_1,\ldots,l_b,\) each greater than or equal to~\(2.\) If~\(b=0\) then  we directly use~(\ref{tj_ind}) to get that
\[\mathbb{E}|w(T_1-\nu)| = \nu^{a} \mathbb{E}|T_1-\nu| \leq \frac{D}{n^{q+1}}\] for some constant~\(D > 0\) and all~\(n\) large.  If~\(b\neq 0,\) then we use the triangle inequality to get
\begin{eqnarray}
\mathbb{E}|w (T_1-\nu)|  &=& \nu^{a}\mathbb{E}\left| T_1\prod_{k=1}^{b} T_{l_k} - \nu\prod_{k=1}^{b} T_{l_k}\right| \nonumber\\
&\leq&  \nu^{a}\mathbb{E}\left| T_1\prod_{k=1}^{b} T_{l_k} - \nu^{b+1}\right| + \nu^{a+1} \mathbb{E}\left|\prod_{k=1}^{b} T_{l_k} - \nu^{b}\right|, \nonumber\\
&\leq& \frac{D}{n^{q+1}} \label{ew_est}
\end{eqnarray}
for some constant~\(D > 0,\) where the first inequality in~(\ref{ew_est}) follows from the triangle inequality and the second inequality in~(\ref{ew_est}) follows from the ``near-independence" estimate~(\ref{et_est2}). The number of terms in the expansion of~\(X\) is at most~\(2^{q}\) and so by the union bound and~(\ref{ew_est}), we get that~\(\mathbb{E}|Q_l| \leq \frac{D \cdot 2^{q}}{n^{q+1}}.\)

Since there are at most~\(n^{q}\) choices for~\(\{i_1,\ldots,i_q\},\) we get from~(\ref{ew_est}) that
\begin{equation}\label{i1_est}
|I_1| \leq \frac{D \cdot (2n)^{q}}{n^{q+1}}.
\end{equation}
Plugging~(\ref{i1_est}) and~(\ref{i2_est}) into~(\ref{etq_split}) we get~\(\mathbb{E}|T-\mu|^{q} \leq n^{\frac{q}{2}} + \frac{D \cdot (2n)^{q}}{n^{q+1}} \leq 2n^{\frac{q}{2}}\) for all~\(n\) large. This proves the second estimate in~(\ref{mu_est}).~\(\qed\)

The following Lemma collects the relevant properties of~\(H  =G(n,p)\) used in the proof of Theorem~\ref{thm_eul_edit}. The complement~\(H^c\) of the graph~\(H\) is defined as follows: For any two vertices~\(u,v\) the edge~\((u,v)\) between~\(u\) and~\(v\) is present in~\(H^c\) if and only if~\((u,v)\) is not present in~\(H.\) A clique in~\(H\) is a complete subgraph of~\(H.\)
\begin{Lemma} \label{t_lem} Suppose~(\ref{p_cond2}) holds and let~\(\epsilon, q > 0\) be constants. There exists an integer~\(n_0 = n_0(q,\epsilon)\) such that if~\(n \geq n_0,\) then the following hold:\\
\((i)\)  If~\(E_{con}(H)\) denotes the event that~\(H\) is connected, then
\begin{equation}\label{f_con}
\mathbb{P}(E_{con}(H)) \geq 1-e^{-\frac{(\log{n})^2}{8}}.
\end{equation}
\((ii)\) If~\(E_{odd}(H)\) denotes the event that the number of odd vertices in the random graph~\(H\) lies between~\(\frac{n}{2}  -n^{\frac{1}{2}+\epsilon}\) and~\(\frac{n}{2} + n^{\frac{1}{2}+\epsilon},\)  then
\begin{equation}\label{e_odd}
\mathbb{P}(E_{odd}(H)) \geq 1-\frac{1}{n^{q\epsilon}}.
\end{equation}
\((iii)\) If~\(E_{good}(H)\) denotes the event that for any two vertices~\(u,v\) there exists a third vertex~\(z\) adjacent to both~\(u\) and~\(v\) in the random graph~\(H,\) then
\begin{equation}\label{e_good}
\mathbb{P}(E_{good}(H) \cap E_{good}(H^c)) \geq 1- 2e^{-(\log{n})^2}.
\end{equation}
\((iv)\) If~\(E_{cliq}(H)\) denotes the event that the largest size of a clique  in~\(H\) is at most~\(2 \sqrt{n} \cdot \log{n},\) then
\begin{equation}\label{e_clique}
\mathbb{P}(E_{cliq}(H) \cap E_{cliq}(H^c)) \geq 1- 2 e^{-(\log{n})^2}.
\end{equation}
\end{Lemma}
Below we prove parts~\((i)\) and~\((ii)\) of Lemma~\ref{t_lem} under the weaker assumption of~(\ref{p_cond}) and prove parts~\((iii)\) and~\((iv)\) assuming~(\ref{p_cond2}). Also throughout, we prove our estimates for the graph~\(H\) first and then use the symmetry of the range of~\(p\) in~(\ref{p_cond2}) and the union bound to obtain the desired estimates in Lemma~\ref{t_lem}.\\
\emph{Proof of Lemma~\ref{t_lem}}\((i)\): We use tree counting arguments to obtain~(\ref{f_con}). Let~\({\cal E}_i\) be the component in~\(H\) containing the vertex~\(i\) and let~\(\#{\cal E}_i\) be the number of vertices in~\({\cal E}_i.\) From estimate~\((3.6)\) in (Ganesan (2018)), we have that
\begin{equation}\label{eir}
\mathbb{P}\left(\#{\cal E}_i  = r\right) \leq {n-1 \choose r-1}\cdot r^{r-2} \cdot p^{r-1}(1-p)^{r(n-r)} \leq (ne)^{r}\cdot e^{-pr(n-r)}
\end{equation}
where the final estimate in~(\ref{eir}) is obtained as follows.  Since~\({n-1 \choose r-1}  = \frac{r}{n} \cdot {n \choose r} \leq {n \choose r},\) we have that~\({n-1 \choose r-1} \cdot r^{r-2} \leq {n \choose r} \cdot r^{r-2}  \leq \left(\frac{ne}{r}\right)^{r} \cdot r^{r-2} \leq (ne)^{r}.\)  Using this estimate together with~\(p^{r-1} \leq 1\) and~\((1-p)^{r(n-r)} \leq e^{-pr(n-r)}\) we get~(\ref{eir}).

For~\(r \leq \frac{n}{2},\) the final term in~(\ref{eir}) is at most~\((ne)^{r} \cdot e^{-\frac{npr}{2}}\) which in turn is bounded above
by~\( \exp\left(-r\left(\frac{np}{2} - \log(ne)\right)\right) \leq \exp\left(-r\left(\frac{(\log{n})^2}{2} - \log(ne)\right)\right)\)
since~\(p \geq \frac{(\log{n})^2}{n}.\) Further,~\(r\left(\frac{(\log{n})^2}{2} - \log(ne)\right) \geq r \cdot \frac{(\log{n})^2}{4} \geq \frac{(\log{n})^2}{4}\) for all~\(n\) large. Thus from~(\ref{eir}) we see that~\(\mathbb{P}\left(\#{\cal E}_i  = r\right) \leq \exp\left(-\frac{(\log{n})^2}{4}\right)\)
for all~\(n\) large. By the union bound, we then get that~\(H\) is disconnected, i.e. that~\(H\) contains a vertex~\(i\) in a component of size at most~\(\frac{n}{2},\) with probability at most~\(n^2 \cdot e^{-\frac{(\log{n})^2}{4}} \leq e^{-\frac{(\log{n})^2}{8}}\) for all~\(n\) large.~\(\qed\)

\emph{Proof of Lemma~\ref{t_lem}\((ii)\)}:  Let~\(\epsilon_1 < \epsilon\) and~\(q_1 > q\) be such that~\(q_1 \epsilon_1 > q \epsilon.\) By Markov's inequality and the second estimate of~(\ref{mu_est}) in Lemma~\ref{vert_lem}, we get for~\(\theta > 0\) that
\[\mathbb{P}\left(|T-\mu| \geq \theta\right) \leq \frac{\mathbb{E}|T-\mu|^{q_1}}{\theta^{q_1}} \leq \frac{D \cdot n^{\frac{q_1}{2}}}{\theta^{q_1}}.\]
Setting~\(\theta = n^{\frac{1}{2} + \epsilon_1}\) and using the mean estimate in~(\ref{mu_est}) we then get~(\ref{e_odd}).~\(\qed\)

\emph{Proof of Lemma~\ref{t_lem}\((iii)\)}: Let~\(u\) and~\(v\) be any two vertices and let~\(z \neq u,v\) be any other vertex. With probability~\(p^2,\) the vertex~\(z\) is adjacent to both~\(u\) and~\(v.\) Therefore with probability~\((1-p^2)^{n-2},\) no vertex in~\(\{1,2,\ldots,n\}\setminus \{u,v\}\) is adjacent to both~\(u\) and~\(v.\) By the union bound, this implies that the event~\(E_{good}(H)\) does not occur with probability
\begin{equation}\label{temp_e_good}
\mathbb{P}(E^c_{comm}(H)) \leq n^2(1-p^2)^{n-2} \leq n^{2} e^{-p^2(n-2)} \leq n^2 e^{-\frac{np^2}{2}} \leq e^{-(\log{n})^2},
\end{equation}
for all~\(n\) large, since~\(p \geq \frac{(\log{n})^2}{\sqrt{n}}\) by~(\ref{p_cond2}). An analogous estimate holds for the complement graph~\(H^c\) and this obtains~(\ref{e_good}).~\(\qed\)

\emph{Proof of Lemma~\ref{t_lem}\((iv)\)}: The probability that there exists a clique of size~\(t\) or more in the random graph~\(H\) is at most
\begin{equation}\label{temp_ab}
{n \choose t} \cdot p^{{t \choose 2}} \leq n^{t} p^{\frac{t(t-1)}{2}} = \left(n\cdot p^{\frac{t-1}{2}}\right)^{t}.
\end{equation}
Using the upper bound for~\(p\) in~(\ref{p_cond2}), we get that
\[n\cdot p^{\frac{t-1}{2}}\leq n\cdot \left(1-\frac{(\log{n})^2}{\sqrt{n}}\right)^{\frac{t-1}{2}} \leq n \cdot \exp\left(-\frac{(\log{n})^2}{\sqrt{n}} \cdot \frac{(t-1)}{2}\right) \leq e^{-(\log{n})^2}\]
for all~\(n\) large, if~\(t = 2\sqrt{n} \cdot \log{n}+1.\) An analogous estimate holds for the complement graph~\(H^c\) and this obtains~(\ref{e_clique}).~\(\qed\)

\renewcommand{\theequation}{\thesection.\arabic{equation}}
\setcounter{equation}{0}
\section{Proof of Theorem~\ref{thm_eul_edit}} \label{thm_proof}
Let~\(H = G(n,p)\) and for~\(\eta > \epsilon,\) let~\(E_{edit}(H)\) be the event that~\(N_{edit}(H)\) lies between~\(\frac{n}{4} -n^{\frac{1}{2}  +\eta}\) and~\(\frac{n}{4}+n^{\frac{1}{2} + \eta}.\)  Similarly define the events~\(E_{ext}(H)\) and~\(E_{red}(H)\) using the quantities~\(N_{ext}(H)\) and~\(N_{red}(H),\) respectively. For any constant~\(\theta > 0\) and for all~\(n\) large, we show below that
\begin{equation}\label{main_eq_h}
\min\left(\mathbb{P}(E_{edit}(H)),\mathbb{P}(E_{ext}(H)),\mathbb{P}(E_{red}(H))\right) \geq 1-\frac{1}{n^{\theta}}
\end{equation}
and this completes the proof of Theorem~\ref{thm_eul_edit}.

In what follows, we discuss edit, extension and reduction numbers in that order. Throughout we use the fact that a graph~\(\Gamma\) is Eulerian if and only if~\(\Gamma\) is connected and each vertex in~\(\Gamma\) has even degree (Theorem~\(1.2.26,\) pp. 27, West (2001)).\\\\
\emph{\underline{Estimate for edits}}: Recalling the events defined in Lemma~\ref{t_lem}, we suppose that the joint event
\begin{equation}\label{e_edit_def}
E_1(H) := E_{con}(H) \cap E_{good}(H) \cap E_{odd}(H)
\end{equation}
occurs. Since~\(E_{con}(H)\) occurs, we know that~\(H\) is connected.

Let~\({\cal T} := \{u_1,\ldots,u_T\}\) be the set of odd degree vertices in~\(H\) so that~\(T\) is an even number. The minimum number of edits needed to convert~\(H\) into an Eulerian graph is~\(\frac{T}{2},\) otherwise there always exists a vertex with odd degree. Using the fact that the event~\(E_{good}(H)\) occurs, we now show that~\(\frac{T}{2}\) edits suffice. Indeed, for~\(1 \leq i \leq \frac{T}{2}\) let~\(e_i\) be the edge with endvertices~\(u_{2i-1}\) and~\(u_{2i}\) in the complete graph~\(K_n\) on~\(n\) vertices. Let~\({\cal F} := \{e_{i_1},\ldots,e_{i_W}\}\) be the set of all edges of~\({\cal E} := \{e_i\}_{1 \leq i \leq \frac{T}{2}}\) present in the graph~\(H\) and consider the graph~\(H_{mod} := \left(H \setminus {\cal F}\right) \bigcup \left({\cal E} \setminus {\cal F}\right)\) obtained by removing the edges in~\({\cal F}\) and adding the edges in~\({\cal E} \setminus {\cal F}\) to the graph~\(H.\) In the graph~\(H_{mod}\) each vertex has even degree. Moreover, because the event~\(E_{good}(H)\) occurs, for each~\(1 \leq i \leq \frac{T}{2}\) there exists a vertex~\(v_i\) that is adjacent to both~\(u_{2i-1}\) and~\(u_{2i}\) in the graph~\(H_{mod}.\) This implies that~\(H_{mod}\) is connected and so~\(H_{mod}\) is Eulerian.

Summarizing we get that~\(N_{edit}(H) = \frac{T}{2}\) and because the event~\(E_{odd}(H)\) also occurs, we have that~\(\frac{T}{2}\) lies between~\(\frac{n}{4}-\frac{n^{\frac{1}{2} + \epsilon}}{2}\) and~\(\frac{n}{4} + \frac{n^{\frac{1}{2} + \epsilon}}{2}.\) Since~\(\eta > \epsilon,\) this implies that the event~\(E_{edit}(H)\) occurs for all~\(n \geq n_0(\eta,\epsilon).\) Finally, for constant~\(\theta > 0\) we see from Lemma~\ref{t_lem} that each event in~(\ref{e_edit_def}) occurs with probability at least~\(1-\frac{1}{n^{\theta}}\) for all~\(n\) large. Therefore by a union bound we get
\begin{equation}\label{e_one_est}
\mathbb{P}(E_{edit}(H)) \geq \mathbb{P}(E_{1}(H)) \geq 1-\frac{3}{n^{\theta}}
\end{equation}
and this completes the estimate for edits.\\\\
\emph{\underline{Estimate for extension}}: Recalling the events defined in Lemma~\ref{t_lem}, suppose that the joint event
\begin{equation}\label{e_ext_def}
E_{2}(H) := E_{con}(H) \cap E_{odd}(H) \cap E_{cliq}(H) \cap E_{good}(H^c)
\end{equation}
occurs. Due to the occurrence of~\(E_{con}(H),\) the graph~\(H\) is connected.

Let~\({\cal T} := \{u_1,\ldots,u_T\}\) be the set of odd degree vertices in~\(H.\) As before at least~\(\frac{T}{2}\) edges must be added to convert~\(H\) into an Eulerian graph. For an upper bound, we argue as follows. Initially all vertices in~\({\cal T}\) are unmarked. If there exists a pair of unmarked vertices~\(u_i,u_j\) in~\({\cal T}\) that are non-adjacent in~\(H,\) we connect these two vertices by an edge and mark~\(u_i\) and~\(u_j.\) We repeat this procedure until we are left with a clique~\(\{v_1,\ldots,v_x\} \subset \{u_1,\ldots,u_T\}\) of unmarked vertices and call the final modified graph as~\(H_{mod}.\) The number~\(x\) is even and the vertices~\(\{v_1,\ldots,v_x\}\) also form a clique in the original graph~\(H.\) Therefore using the fact that~\(E_{cliq}(H)\) occurs, we get that~\(x \leq 2 \sqrt{n} \cdot \log{n}.\) We now divide the vertices in~\(\{v_1,\ldots,v_x\}\) into~\(\frac{x}{2}\) pairs~\(\{v_{2i-1},v_{2i}\}_{1 \leq i \leq \frac{x}{2}}.\) Since~\(E_{good}(H^c)\) occurs, each pair~\(\{v_{2i-1},v_{2i}\}\) contains a common neighbour~\(z_i\) in the complement graph~\(H^c.\)  Adding the edges~\(\{(v_{2i-1},z_i),(v_{2i},z_i)\}_{1 \leq i \leq \frac{x}{2}}\) to the graph~\(H_{mod}\) we then get a new connected graph~\(G_{mod}\) where each vertex has even degree.

The graph~\(G_{mod}\) is Eulerian and is obtained after adding at most~\(\frac{T}{2} + x\) edges to the original graph~\(G.\) Summarizing we have that if the event~\(E_{ext}\) defined in~(\ref{e_ext_def}) occurs, then~\(G\) can be converted into an Eulerian graph by adding at most~\(\frac{n}{4} + \frac{n^{\frac{1}{2}+\epsilon}}{2} + 2\sqrt{n} \cdot \log{n}\) edges. Recalling that~\(\eta > 0\) is larger than~\(\epsilon,\) we then get that the event~\(E_{ext}(H)\) occurs for all~\(n \geq n_0(\eta, \epsilon).\) As before, we use Lemma~\ref{t_lem} to estimate the probabilities of the events in~(\ref{e_ext_def}) to get for~\(\theta > 0\) and all~\(n\) large that
\begin{equation}\label{e_two_est}
\mathbb{P}(E_{ext}(H)) \geq \mathbb{P}(E_{2}(H)) \geq 1-\frac{4}{n^{\theta}}.
\end{equation}
This completes the estimate for extensions.\\\\
\emph{\underline{Estimate for reduction}}: Suppose that the joint event
\begin{equation}\label{e_red_def}
E_{3}(H) := E_{odd}(H) \cap E_{con}(H) \cap E_{good}(H) \cap E_{cliq}(H^c)
\end{equation}
occurs. Due to the occurrence of~\(E_{cliq}(H^c),\) the maximum size of a clique in~\(H^c\) is at most~\(2 \sqrt{n} \cdot \log{n}.\) This implies that the maximum size of an independent set in~\(H\) is at most~\(2 \sqrt{n} \cdot \log{n},\) where we recall that a set of vertices~\({\cal V}\) is said to be an independent set in~\(H\) if no two vertices of~\({\cal V}\) are adjacent in~\(H.\)

As before, let~\({\cal T} := \{u_1,\ldots,u_T\}\) be the set of odd degree vertices in~\(H\) so that at least~\(\frac{T}{2}\) edges must be removed to convert~\(H\) into an Eulerian graph. For obtaining an upper bound, we proceed similar to the extension case. Initially all vertices in~\({\cal T}\) are unmarked. If there exists a unmarked pair of vertices~\(u_i,u_j\) in~\({\cal T}\) that are adjacent in~\(H,\) we remove the edge~\((u_i,u_j)\) to get a new graph~\(F_{temp}\) and mark~\(u_i\) and~\(u_j.\) Because~\(E_{good}(H)\) occurs, there exists a common neighbour~\(z_{ij}\) of the vertices~\(u_i\) and~\(u_j\) in the graph~\(H\) and so~\(F_{temp}\) remains connected. We now repeat this procedure until we are left with an independent set~\(\{w_1,\ldots,w_y\} \subset \{u_1,\ldots,u_T\}\) of unmarked vertices and call the final modified graph as~\(F_{mod}.\) Arguing as above, the graph~\(F_{mod}\) still remains connected.

As in the extension case, we see that the number~\(y\) is even and because~\(E_{cliq}(H^c)\) occurs, we must have that~\(y \leq 2 \sqrt{n} \cdot \log{n}.\) We now divide the vertices in~\(\{w_1,\ldots,w_y\}\) into~\(\frac{y}{2}\) pairs~\(\{w_{2i-1},w_{2i}\}_{1 \leq i \leq \frac{y}{2}}\) and iteratively remove edges attached to each pair. Using the  fact that~\(E_{good}(H)\) occurs, we see that each pair~\(\{w_{2i-1},w_{2i}\}_{1 \leq i \leq \frac{y}{2}}\) contains a common neighbour~\(s_{2i-1,2i}\) in the graph~\(H.\)  By construction, the vertex~\(s_{2i-1,2i}\) is adjacent to~\(w_{2i-1}\) and~\(w_{2i}\) in the graph~\(F_{mod}\) as well.

First we remove the edges~\((w_1,s_{12})\) and~\((w_2,s_{12})\) from~\(F_{mod}\) to obtain a graph~\(F_{new}.\) Again using the fact that~\(E_{good}(H)\) occurs, we know that the vertices~\(s_{12}\) and~\(w_1\) share a common neighbour~\(b_1\) in the graph~\(H,\) which does not belong to the independent set~\(\{w_1,\ldots,w_y\}.\) Therefore the edges~\((b_1,s_{12})\) and~\((b_1,w_1)\) are present in the graph~\(F_{new}\) as well. Similarly~\(s_{12}\) and~\(w_2\) share a common neighbour~\(b_2\) in~\(F_{new}\) and so~\(F_{new}\) remains connected.

We now proceed iteratively and remove the edges\\\(\{(w_{2i-1},s_{2i-1,2i}),(w_{2i},s_{2i-1,2i})\}_{1 \leq i \leq \frac{y}{2}}\) from the graph~\(F_{mod}\) to get a new graph\\\(G_{mod}\) containing only even degree vertices. Arguing as in the previous paragraph, we see that the graph~\(G_{mod}\) still remains connected.  We then argue as in the extension case (see discussion preceding~(\ref{e_two_est})) to get that~\(E_{red}(H)\) occurs for all~\(n \geq n_0(\eta,\epsilon)\) and that
\begin{equation}\label{e_three_est}
\mathbb{P}(E_{red}(H)) \geq \mathbb{P}(E_{3}(H)) \geq 1-\frac{4}{n^{\theta}}.
\end{equation}
This completes the estimate for reductions.

From~(\ref{e_one_est}),~(\ref{e_two_est}) and~(\ref{e_three_est}), we get~(\ref{main_eq_h}).~\(\qed\)


%






\subsection*{Acknowledgement}
I thank Professors V. Raman, C. R. Subramanian and the referees for crucial comments that led to an improvement of the paper. I also thank IMSc for my fellowships.

\bibliographystyle{plain}

\end{document}